# Doubling properties for second order parabolic equations

*In memory of Eugene Fabes*

By Mikhail V. Safonov and Yu Yuan*

## Abstract

We prove the doubling property of $L$-caloric measure corresponding to the second order parabolic equation in the whole space and in Lipschitz domains. For parabolic equations in the divergence form, a weaker form of the doubling property follows easily from a recent result, the backward Harnack inequality, and known estimates of Green's function. Our method works for both the divergence and nondivergence cases. Moreover, the backward Harnack inequality and estimates of Green's function are not needed in the course of proof.

## 1. Introduction

A measure is said to be doubling, or to satisfy the doubling condition, if for any pair of concentric balls with radii $r$ and $2r$, their measures are comparable. The doubling property for certain measures is the starting point in harmonic analysis, for example, in deriving the weak-type $(1,1)$ estimate and $L^p$ inequalities for the maximal operator. The doubling property for parabolic and elliptic equations is also essential in extending the classical Fatou theorem. It is known that the doubling property for elliptic equations holds true. See [Ke] and other references below for further details and applications.

The purpose of this paper is to establish the doubling property of $L$-caloric measure corresponding to the second order parabolic equation $Lu = 0$ in a cylinder $Q = \Omega \times (0, \infty)$ with $\Omega$ being $\mathbb{R}^n$ or a Lipschitz domain in $\mathbb{R}^n$ with Lipschitz constants $m$, $r_0$ (see Assumptions in Section 2). We consider both the divergence (D) and nondivergence (ND) operators $L$:

*\**Key words and phrases*. Doubling property, $L$-caloric measure, Fatou theorem.
Both authors are partially supported by NSF Grant No. DMS-9623287.

$$\text{(D)} \qquad Lu = \sum_{i,j=1}^{n} D_i\left(a_{ij}(x,t) D_j u(x,t)\right) - D_t u(x,t),$$

$$\text{(ND)} \qquad Lu = \sum_{i,j=1}^{n} a_{ij}(x,t) D_{ij} u(x,t) - D_t u(x,t),$$

where $D_j u = \partial u/\partial x_j$, $D_{ij} u = D_i D_j u$, $D_t u = \partial u/\partial t$. We assume that the coefficients $a_{ij} = a_{ij}(x,t) \in C^{\infty}(\overline{Q})$, and also for all $X = (x,t) \in Q$, $\xi = (\xi_1, \cdots, \xi_n) \in \mathbb{R}^n$,

$$\nu|\xi|^2 \leq \sum_{i,j=1}^{n} a_{ij}(X)\xi_i\xi_j, \qquad \max_{i,j}|a_{ij}(X)| \leq \nu^{-1},$$

with a constant $\nu \in (0,1]$. The assumption $a_{ij} \in C^{\infty}$ is qualitative, in the sense that none of our estimates depend on the smoothness of $a_{ij}$. By the standard approximation technique, all of our results are extended to the parabolic equations with measurable $a_{ij}$ in the divergence case, and with continuous $a_{ij}$ in the nondivergence case.

For any continuous function $\varphi$ on the parabolic boundary $\partial_p Q$ (which is defined in Section 2), the bounded solution of the boundary value problem

$$Lu = 0 \quad \text{in} \quad Q, \qquad u = \varphi \quad \text{on} \quad \partial_p Q$$

is represented by means of the $L$-caloric measure $\omega^X = \omega_Q^X$ on $\partial_p Q$ as follows:

$$u(X) = u(x,t) = \int_{\partial_p Q} \varphi(Y) d\omega^X(Y)$$

(see [A], [FGS], [G], [LSU]).

The following two theorems contain our main results concerning the doubling properties of the $L$-caloric measure in $Q = \Omega \times (0, \infty)$ where $\Omega$ is the whole space $\mathbb{R}^n$ or a Lipschitz domain in $\mathbb{R}^n$.

THEOREM 1.1. *Let the coefficients of operator $L$ be defined on $Q = \mathbb{R}^n \times (0, \infty)$, and let a constant $K \geq 1$ be fixed. Then for all $r > 0$ and $X = (x,t) \in Q$ with $|x| \leq K\sqrt{t}$,*

$$(1.1) \qquad \omega^X(\triangle_{2r}) \leq N\omega^X(\triangle_r),$$

*where $\triangle_r = B_r(0) \times \{0\} \subset \mathbb{R}^n \times \{0\} = \partial_p Q$, and the constant $N = N(n, \nu, K)$.*

THEOREM 1.2. *Let the coefficients of operator $L$ be defined on $Q = \Omega \times (0, \infty)$, where $\Omega$ is a Lipschitz domain in $\mathbb{R}^n$ with Lipschitz constants $m, r_0$. Let $Y = (y,s) \in \partial_p Q$ and constants $K \geq 1, \lambda \geq 1$ be given. Then for all $r \in (0, \lambda r_0/4]$ and $X = (x,t) \in Q$ satisfying*

$$(1.2) \qquad |x - y| \leq K\sqrt{t-s}, \qquad 4r \leq \sqrt{t-s} \leq \lambda r_0,$$

*the estimate* (1.1) *holds, where*

$$\triangle_r = \triangle_r(Y) = \{Z = (z, \tau) \in \partial_p Q : |z - y| < r, |\tau - s| < r^2\},$$

*and the constant* $N = N(n, \nu, m, \lambda, K)$.

*Remark* 1.1. For a class of unbounded Lipschitz domains with $r_0 = \infty$, e.g. half space, Theorem 1.2 holds true with $\lambda = 1$. For bounded domains $\Omega$, we have the estimate (1.1) with $N$ corresponding to $\lambda r_0 = \operatorname{diam} \Omega / K = R$, for all $r \in (0, R/4]$ and $X \in \Omega \times \{s + R^2\}$. By the comparison principle, this estimate is extended to all $X \in \Omega \times (s + R^2, \infty)$, with the same constant $N$. Note that $r$ cannot be arbitrarily large for bounded domains, for easy examples show that $N \to \infty$ as $r/r_0 \to \infty$. We can also derive the conclusion of Theorem 1.2 for bounded domain $\Omega$, and $r$, $X = (x, t)$ satisfying

(1.3) $\quad t - s \geq \delta^2 > 0, \ 0 < r \leq \dfrac{1}{4} \min(r_0, \delta) \quad (\delta = \operatorname{const} > 0),$

instead of (1.2), because (1.2) follows from (1.3) with $K = \operatorname{diam} \Omega / \delta$.

In the divergence case, Theorem 1.1 is an immediate consequence of Aronson's estimate (see [A]) for the fundamental solution. Theorem 1.2 was proved in [FGS, Th. 2.4] in the divergence, time-independent case for bounded domains $\Omega$, with the conditions (1.3) instead of (1.2), by means of the estimates for the corresponding Green's function. It was also outlined in [FGS] that in this case, Theorem 1.2 is actually equivalent to the backward Harnack inequality which is formulated in Theorem 2.3 below and was recently proved in [FS]. Thus in the divergence case one can derive Theorem 1.1 and a "weak" form of Theorem 1.2 for bounded domain $\Omega$, with (1.3) instead of (1.2), from the results in [A], [FGS] and [FS]. We demonstrate this approach here in Section 3, after some preparations in Section 2 where we introduce notation and collect together known results. Along the same lines, the doubling property was proved for the divergence equations with singular drift terms in [HL].

In the nondivergence case, the appropriate estimates for the fundamental solution and Green's function fail (see [FK], [S]), and the backward Harnack inequality does not help in the proof of the doubling property. In Section 4, we present an alternative approach which works for both the divergence and nondivergence cases in full generality and which does not use the backward Harnack inequality. Moreover, even the usual Harnack inequality (Theorem 2.2) is not needed in the proof of Theorem 1.1 in the nondivergence case. Instead we could apply the comparison principle (Theorem 2.1) in combination with some simple barrier functions. However, we prefer using the Harnack inequality, which simplifies the proofs and makes it possible to treat simultaneously both the divergence and nondivergence cases.

One of applications of the doubling property is the *Fatou theorem* (Theorem 2.14 in [FGS]) which states that any positive solution of $Lu = 0$ in $Q = \Omega \times (0, \infty)$ has finite nontangential limits at almost every (with respect to the $L$-caloric measure $\omega$) point $Y \in \partial\Omega \times (0, \infty)$. Following the framework of the paper [FGS], where the Fatou theorem was proved in the divergence, time-independent case, and also the papers [B], [CFMS], [FGMS] dealing with the elliptic equations, one can extend the Fatou theorem to the general divergence and nondivergence cases. For other applications, we refer the reader to [JK], [Ke], [ACS], [HL], and the references therein.

Throughout this paper, $N$ will denote various positive constants depending only on the original quantities.

## 2. Notation and known results

For an arbitrary domain $V \subset \mathbb{R}^{n+1}$, we define its *parabolic boundary* $\partial_p V$ as the set of all the points $X \in \partial V$ such that there is a continuous curve lying in $V \cup \{X\}$ with initial point $X$, along which $t$ is nondecreasing. In particular, for $Q = \Omega \times (0, T)$ we have

$$\partial_p Q = \partial_x Q \cup \partial_t Q, \text{ where } \partial_x Q = \partial\Omega \times (0, T), \ \partial_t Q = \overline{\Omega} \times \{0\}.$$

For $\delta = \text{const} > 0$, $\Omega \subset \mathbb{R}^n$, $Q = \Omega \times (0, T)$, we set

(2.1) $\qquad \Omega^\delta = \{x \in \Omega : \text{dist}(x, \partial\Omega) > \delta\}, \qquad Q^\delta = \Omega^\delta \times (\delta^2, T).$

For $y \in \mathbb{R}^n$, $r > 0$, $\delta > 0$, and a domain $\Omega \subset \mathbb{R}^n$,

$$\begin{aligned} B_r(y) &= \{x \in \mathbb{R}^n : |x - y| < r\}, \qquad \Omega_r = \Omega_r(y) = \Omega \cap B_r(y), \\ \Omega_r^\delta &= \Omega_r^\delta(y) = \Omega^\delta \cap B_r(y) = \{x \in \Omega : \text{dist}(x, \partial\Omega) > \delta, \ |x - y| < r\}. \end{aligned}$$

For $Y = (y, s) \in \mathbb{R}^{n+1}$, $r > 0$, and a cylinder $Q = \Omega \times (0, \infty)$,

$$Q_r = Q_r(Y) = Q \cap C_r(Y), \qquad \triangle_r = \triangle_r(Y) = (\partial_p Q) \cap C_r(Y),$$

where $C_r(Y) = B_r(y) \times (s - r^2, s + r^2)$ is a "standard" cylinder. In particular, $\triangle_r$ in Theorem 1.1 corresponds to $Y = 0 \in \mathbb{R}^{n+1}$. We will also use more general cylinders

$$C_{R,r}(Y) = B_R(y) \times (s - r^2, s + r^2), \quad Q_{R,r}(Y) = Q \cap C_{R,r}(Y).$$

The following *comparison principle* is well-known.

THEOREM 2.1. *Let $V$ be a bounded domain in $\mathbb{R}^{n+1}$, and let functions $u, v \in C^2(V) \cap C(\overline{V})$ satisfy $Lu \leq Lv$ in $V$, $u \geq v$ on $\partial_p V$. Then $u \geq v$ on $\overline{V}$.*

THEOREM 2.2 (Harnack inequality). *Let $u$ be a nonnegative solution of $Lu = 0$ in a bounded cylinder $Q = \Omega \times (0, T)$, $\delta = \text{const} > 0$ be such that* $\Omega^\delta$ is a connected set, $\operatorname{diam}\Omega/\delta \leq \lambda$, $T/\delta^2 \leq \lambda = const < \infty$. Then for all $X = (x,t), Y = (y,s) \in Q^\delta$ satisfying $t - s \geq \delta^2$,

$$u(Y) \leq Nu(X), \tag{2.2}$$

where the constant $N = N(n,\nu,\lambda)$.

The above theorem in the divergence case was proved in [M1]; see also [A], [FSt]. In the nondivergence case it was proved in [KS]; see also [Kr, Ch. 4].

From now on we assume that the domain $\Omega \subset \mathbb{R}^n$ satisfies the following *Lipschitz condition* with some positive constants $r_0$, $m$.

*Assumptions*: For each $y \in \partial\Omega$, there is an orthonormal coordinate system (centered at $y$) such that

$$\Omega \cap B_{r_0}(y) = \{x = (x', x_n) : x' \in \mathbb{R}^{n-1}, x_n > \varphi(x'), |x| < r_0\},$$

where $\|\nabla\varphi\|_{\mathfrak{L}^\infty} \leq m$.

Under these assumptions, there exists a constant $\mu = \mu(m) > 0$ such that for arbitrary $y \in \partial\Omega$ and $0 < r < r_0$, we have

$$B_{2\mu r}(y^{(1)}) \subset B_r(y) \cap \Omega, \qquad B_{2\mu r}(y^{(2)}) \subset B_r(y) \setminus \Omega \tag{2.3}$$

for some $y^{(1)}, y^{(2)} \in \mathbb{R}^n$.

The next result is recent, the backward Harnack inequality. It is contained in [FS, Th. 1.4] and [FSY, Th. 3.7]. We formulate it here in an equivalent form.

THEOREM 2.3. *Let $u$ be a nonnegative solution of $Lu = 0$ in $Q = \Omega \times (0,\infty)$ which continuously vanishes on $\partial_x Q = \partial\Omega \times (0,\infty)$, and let constants $\delta > 0$, $\mu > 0$, $\operatorname{diam}\Omega/r_0 \leq \lambda$, $\operatorname{diam}\Omega/\delta \leq \lambda$. Then for $X = (x,t) \in Q$ and $r$ satisfying*

$$\operatorname{dist}(x, \partial\Omega) > \mu r, \quad t > \delta^2, \quad 0 < r \leq \frac{1}{2}\min(r_0, \delta),$$

*we have*

$$u(x, t+r^2) \leq Nu(x, t-r^2),$$

*where the constant $N = N(n,\nu,m,\lambda,\mu)$.*

The following theorem helps to control the quotient of two solutions near $\partial_x Q$. It originates in [FGS] and [FSY].

THEOREM 2.4. *Let $u$ and $v$ be two positive solutions of $Lu = 0$ in $Q = \Omega \times (0,\infty)$ which continuously vanish on $(\partial_x Q) \cap C_{(K+2)r, 2r}(Y)$, where $Y = (y,s) \in \overline{Q}$, $K \geq 1$, $0 < r \leq r_0/4$, and $s \geq 5r^2$. Then*

$$\sup_{Q_{Kr,r}(Y)} \frac{v}{u} \leq N \frac{\inf_{\Omega_{2r}^+} v}{\sup_{\Omega_{2r}^-} u}, \tag{2.4}$$

*where*

$$\Omega_\rho^\pm = \Omega_\rho^\pm(Y) = \Omega_\rho^{\mu\rho}(y) \times \{s \pm \rho^2\} = (\Omega^{\mu\rho} \cap B_\rho(y)) \times \{s \pm \rho^2\}$$

*for $0 < \rho \leq r_0/2$, and the constant $N = N(n, \nu, m, K)$. Here $\mu = \mu(m) > 0$ is the constant in (2.3).*

*Proof.* For arbitrary $X = (x, t) \in Q_{Kr,r}(Y)$, we first consider the case dist $(x, \partial\Omega) < r$. Then $X \in Q_r(X_0)$ for some $X_0 = (x_0, t) \in \partial_x Q$. By our assumptions, $u = 0$ on $(\partial_x Q) \cap C_r(Y)$. From [FGS, Th. 1.6] in the divergence case, and [FSY, Th. 4.3] in the nondivergence case, it follows that there exists a constant $\varepsilon = \varepsilon(n, \nu, m) > 0$ and points $X_r^\pm = (x_r, t \pm r^2)$ with $x_r \in \Omega_r^{\varepsilon r}(x_0)$ such that

$$\sup_{Q_{\varepsilon r}(X_0)} \frac{v}{u} \leq N(n, \nu, m) \frac{v(X_r^+)}{u(X_r^-)}.$$

Further, by the Harnack inequality,

$$v(X_r^+) \leq N \inf_{\Omega_{2r}^+} v, \quad \sup_{\Omega_{2r}^-} u \leq Nu(X_r^-).$$

These estimates yield (2.4).

If dist $(x, \partial\Omega) \geq r$, we can apply the Harnack inequality directly, which implies

$$v(X) \leq N \inf_{\Omega_{2r}^+} v, \quad \sup_{\Omega_{2r}^-} u \leq Nu(X),$$

and we also have (2.4). Hence the estimate (2.4) holds for all $X = (x, t) \in Q_{Kr,r}(Y)$. □

## 3. The divergence case

In this section, we sketch the proofs of Theorem 1.1 and a special case of Theorem 1.2, with (1.3) instead of (1.2), in the divergence case only. Our approach here follows [FGS] and is based on the Gaussian estimates for the fundamental solution and Green's function. Such estimates fail in the nondivergence case. In Section 4, we prove these theorems again by a more general method, which works in both the divergence and nondivergence cases simultaneously and does not need the additional restriction (1.3).

3.1. *Proof of Theorem* 1.1. By the substitutions, $x \longrightarrow \lambda x$, $t \longrightarrow \lambda^2 t$, where $\lambda = $ const, we reduce the proof to the case $X = (x, 1)$. Then

$$\omega^X(\triangle_r) = \int_{B_r(0)} \Gamma(x, 1; y, 0)\, dy,$$

where $\Gamma$ is the fundamental solution corresponding to the divergence parabolic operator $L$. By Aronson's estimate ([A, Th. 7]),

$$(3.1) \quad \frac{1}{N}\exp\left(-N|x-y|^2\right) \leq \Gamma(x,1;y,0) \leq N\exp\left(-\frac{|x-y|^2}{N}\right),$$

where $N = N(n,\nu)$. If $r \leq 1$, then from our assumption $|x| \leq K$ it follows that the above exponents lie between two positive constants for $y \in B_{2r}(0)$. Hence

$$\omega^X(\triangle_{2r}) \leq Nr^n, \qquad r^n \leq N\omega^X(\triangle_r),$$

and we get the desired estimate (1.1). If $r > 1$, then

$$\omega^X(\triangle_{2r}) \leq 1 \leq N\omega^X(\triangle_1) \leq N\omega^X(\triangle_r);$$

i.e. we also have (1.1). □

3.2. *Proof of Theorem* 1.2. We give only the outline of the proof of this theorem with conditions (1.3) and bounded $\Omega$, because it is quite similar to the proof of Theorem 2.4 in [FGS] where the doubling property is stated in the time-independent case.

First of all, using scaling, we reduce the proof to the case $\operatorname{diam}\Omega = 1$. We take $Y = (y,s) \in \partial_x Q = \partial\Omega \times (0,\infty)$. The case $Y = (y,s) \in \partial_t Q = \overline{\Omega} \times \{0\}$ can be treated with the same technical adjustments as in [FGS]. Moreover, by the comparison principle, it suffices to consider $X = (x,t)$ with $t - s \leq 1$, so that we can restrict ourselves to a bounded cylinder $Q = \Omega \times (0,T) \subset \mathbb{R}^{n+1}$.

We can always assume that $r_0 \leq \delta$ and the coefficients $a_{ij}(x,t)$ of $L$ are extended for $t < 0$, and hence Green's function $G(X;Y) = G(x,t;y,s)$ is well-defined in the cylinder $\Omega \times (-T,T)$. Following the lines of the proof of Theorem 1.4 in [FGS], we have

$$N^{-1}\rho^n G(X;Y_\rho^+) \leq \omega^X(\triangle_\rho(Y)) \leq N\rho^n G(X;Y_\rho^-)$$

for $0 < \rho \leq r_0/2$, where $Y_\rho^\pm = (y_\rho, s \pm \rho^2)$, $y_\rho \in \Omega_\rho^{\mu\rho}(y)$, and the constant $\mu = \mu(m) \in (0,1)$.

Further, for fixed $X = (x,t)$, the function $G(X;Y) = G(x,t;y,s)$ is a solution of the parabolic equation with $L$-adjoint operator $L^* = \sum D_j(a_{ij}D_i) + D_s$, where $D_i = \partial/\partial y_i$. Substituting $s$ by $-s$, one can reduce this operator to the same form as the divergence operator $L$ in (D). Applying Theorems 2.3 and 2.2, we obtain

$$N^{-1}G(X;Y_{2r}^-) \leq G(X;Y_{2r}^+) \leq NG(X;Y_r^+)$$

for $0 < r \leq r_0/4$. These estimates give us

$$\begin{aligned}\omega^X(\triangle_{2r}(Y)) &\leq N_1 r^n G(X;Y_{2r}^-)\\ &\leq N_2 r^n G(X;Y_r^+) \leq N\omega^X(\triangle_r(Y));\end{aligned}$$

i.e. we have the desired estimate (1.1). □

## 4. The general case

It is now our intention to prove Theorem 1.1 and Theorem 1.2 in both the divergence and nondivergence cases corresponding to operators $L$ in (D) and (ND). The operators $L$ are defined in the cylinder $Q = \Omega \times (0, \infty)$, where $\Omega$ is the whole space $\mathbb{R}^n$ in Theorem 1.1, and is a Lipschitz domain in $\mathbb{R}^n$ with positive constants $r_0, m$ in Theorem 1.2.

4.1. *Auxiliary results.* In the following lemmas, statements (a) and (b) correspond respectively to the cases: (a) $\Omega = \mathbb{R}^n$, and (b) $\Omega$ is a Lipschitz domain in $\mathbb{R}^n$.

The next lemma contains a standard estimate for $L$-caloric measures. Usually such estimates are proved by different means in the divergence and nondivergence cases. Here we give a proof which is based only on the Harnack inequality and is valid for both these two cases.

LEMMA 4.1.   (a) *Let* $Q = \mathbb{R}^n \times (0, \infty)$, $Y = (y, 0) \in \partial_p Q$, *and* $r > 0$. *There exist a positive constant* $N = N(n, \nu)$ *such that the $L$-caloric measure* $\omega^X$ *satisfies*

$$(4.1) \qquad \inf_{Q_r(Y)} \omega^X(\triangle_{2r}) \geq \frac{1}{N},$$

*where* $Q_r(Y) = Q \cap C_r(Y) = B_r(y) \times (0, r^2)$, $\triangle_{2r} = \triangle_{2r}(Y) = (\partial_p Q) \cap C_{2r}(Y)$.

(b) *Let* $Q = \Omega \times (0, \infty)$, $Y = (y, s) \in \partial_p Q$, *and* $0 < r \leq r_0/2$. *Then the estimate* (4.1) *holds with* $N = N(n, \nu, m)$.

*Proof.* (a) Let $\omega_C^X$ denote the $L$-caloric measure for $C = C_{2r}(Y)$, and $v(X) = \omega_C^X(\partial_p C \cap \{t \leq 0\})$. Then automatically $v \equiv 1$ on $C \cap \{t \leq 0\}$. By the comparison principle (Theorem 2.1), we get

$$\omega^X(\triangle_{2r}) \geq v(X) \quad \text{in} \quad Q_{2r}(Y).$$

Using the Harnack inequality (Theorem 2.2) applied to $v$ in $C$, we get the desired estimate (4.1):

$$\inf_{Q_r(Y)} \omega^X(\triangle_{2r}) \geq \inf_{Q_r(Y)} v \geq \frac{1}{N} v(y, -r^2) = \frac{1}{N}.$$

(b) The above proof of the statement (a) is valid for $Y = (y, 0) \in \partial_t Q = \overline{\Omega} \times \{0\}$ without any modifications, so that it remains to consider the case $Y = (y, s) \in \partial_x Q = \partial\Omega \times (0, \infty)$. By the properties (2.3), there exists a cylinder

$$C' = B_{\mu r}(z) \times (s - 4r^2, s + 4r^2) \subset C \setminus Q = C_{2r}(Y) \setminus Q,$$

where the constant $\mu = \mu(m) > 0$. Let $\triangle'_{\mu r} = B_{\mu r}(r) \times \{s - 4r^2\}$ denote the bottom of this cylinder. Using the comparison principle twice, in $C$ and in $C'$,

we have
$$\omega^X(\triangle_{2r}) \geq v(X) = \omega_C^X(\triangle'_{\mu r}) \quad \text{in} \quad Q_{2r}(Y),$$
$$v(X) \geq v'(X) = \omega_{C'}^X(\triangle'_{\mu r}) \quad \text{in} \quad C'.$$

By the Harnack inequality applied to $v$ in $C$, we get
$$\inf_{Q_r(Y)} \omega^X(\triangle_{2r}) \geq \inf_{Q_r(Y)} v \geq \frac{1}{N} v(z, s - 2r^2) \geq \frac{1}{N} v'(z, s - 2r^2).$$

One can extend $v'$ from $C'$ to a longer cylinder
$$C'' = B_{\mu r}(z) \times (s - 5r^2, s + 4r^2)$$
by the formula
$$v'(X) = \omega_{C''}^X(\partial_p C'' \cap \{t \leq s - 4r^2\}),$$
so that $v' \equiv 1$ on $C'' \cap \{t \leq s - 4r^2\}$. Applying the Harnack inequality to $v'$ in $C''$, we obtain
$$v'(z, s - 2r^2) \geq \frac{1}{N} v'(z, s - 4r^2) = \frac{1}{N}.$$
This inequality together with the previous one yields (4.1). □

COROLLARY 4.2. *Under the assumptions of Lemma 4.1, let $u$ be a solution of $Lu = 0$ which continuously vanishes on $\triangle_{2r} = (\partial_p Q) \cap C_{2r}(Y)$. Then the positive and negative parts of $u$, $u^\pm = \max(\pm u, 0)$, satisfy*

(4.2) $$\sup_{Q_r(Y)} (u^\pm) \leq \theta \sup_{Q_{2r}(Y)} (u^\pm)$$

*with a constant $\theta = \theta(n, \nu, m) \in (0, 1)$. In the case $\Omega = \mathbb{R}^n$, the constant $\theta$ does not depend on $m$.*

*Proof.* For arbitrary $X \in Q_r = Q_r(Y)$,
$$u(X) = \int_{\partial Q_{2r}} u \, d\omega^X = \int_{(\partial Q_{2r}) \setminus \triangle_{2r}} u \, d\omega^X.$$

Hence
$$u^\pm(X) = \max(\pm u(X), 0) \leq \int_{(\partial Q_{2r}) \setminus \triangle_{2r}} (u^\pm) \, d\omega^X$$
$$\leq \omega^X((\partial Q_{2r}) \setminus \triangle_{2r}) \cdot \sup_{Q_{2r}}(u^\pm) = (1 - \omega^X(\triangle_{2r})) \cdot \sup_{Q_{2r}}(u^\pm),$$

and hence (4.1) implies (4.2) with $\theta = 1 - N^{-1} < 1$. □

LEMMA 4.3. (a) *Let $Q = \mathbb{R}^n \times (0, \infty)$, $Y = (y, 0) \in \partial_p Q$, and let $u$ be a solution of $Lu = 0$ in $Q$ such that*

(4.3) $$u \geq 0 \quad \text{on} \quad U_R = \{(x, t) : |x - y| \leq K\sqrt{t} \leq KR\},$$

where $K \geq 8$ and $R > 0$ are given constants. Then the function

$$f_1(\rho) = \inf_{B_\rho^+} u, \quad \text{where} \quad B_\rho^+ = B_\rho(0) \times \{\rho^2\},$$

satisfies

(4.4) $\quad f_1(\rho) \geq \left(\dfrac{\rho_0}{\rho}\right)^{\gamma_1} \inf_{\rho_0 \leq r \leq 2\rho_0} f_1(r) \quad \text{for} \quad 0 < 2\rho_0 \leq \rho \leq R,$

with a constant $\gamma_1 = \gamma_1(n, \nu) > 0$.

(b) Let $Q = \Omega \times (0, \infty)$, $Y = (y, s) \in \partial_p Q$, and let $u$ be a solution of $Lu = 0$ in $Q$ such that

(4.5) $\quad u \geq 0 \text{ on } U_R' = Q \cap \left\{|x - y| \leq K\sqrt{t - s},\ \rho_0 \leq \sqrt{t - s} \leq R\right\}$

with constants $K \geq 8$, $0 < 2\rho_0 \leq R \leq \lambda r_0$. Then the function

$$f_1(\rho) = \inf_{\Omega_\rho^+} u, \quad \text{where} \quad \Omega_\rho^+ = \Omega_\rho^{\mu\rho'} \times \{s + \rho^2\},\ \rho' = \min(\rho, r_0),$$

and $\mu = \mu(m) > 0$ is the constant in (2.3), satisfies (4.4) with a constant $\gamma_1 = \gamma_1(n, \nu, m, \lambda) > 0$.

*Proof.* (a). For arbitrary $\rho \in (2\rho_0, R]$, the sets $B_{\rho/2}^+$, $B_\rho^+$ lie in the closure $\overline{C}$ of the cylinder

$$C = B_{(1+\varepsilon)\rho}(y) \times ((1 - \varepsilon)\rho^2/4, \rho^2) \subset U_R$$

with a small absolute constant $\varepsilon > 0$, and stay away from its parabolic boundary $\partial_p C$. Hence we can apply the Harnack inequality which gives us

$$f_1(\rho) = \inf_{B_\rho^+} u \geq 2^{-\gamma_1} \sup_{B_{\rho/2}^+} u \geq 2^{-\gamma_1} f_1(\rho/2)$$

with $\gamma_1 = \gamma_1(n, \nu) > 0$. Iterating this inequality $k$ times, where $k$ satisfies $2\rho_0 > 2^{-k}\rho \geq \rho_0$, we see that

$$f_1(\rho) \geq 2^{-k\gamma_1} f_1(2^{-k}\rho) \geq \left(\dfrac{\rho_0}{\rho}\right)^{\gamma_1} f_1(2^{-k}\rho)$$

which implies (4.4). Thus statement (a) is proved.

The proof of (b) is essentially the same, only $B^+(\rho)$ should be replaced by $\Omega_\rho^+$, and the cylinder $C$ by the cylinder

$$\Omega_{(1+\varepsilon)\rho}^{\mu\rho'/2}(y) \times ((1 - \varepsilon)\rho^2/4, \rho^2) \subset U_R'.$$

This completes the proof of Lemma 4.3. $\square$

LEMMA 4.4. *Let $Q = \Omega \times (0, \infty)$, where $\Omega$ is either $\mathbb{R}^n$ or a Lipschitz domain in $\mathbb{R}^n$, and $Y = (y, s) \in \partial_p Q$. Let $u$ be a solution of $Lu = 0$ in $Q$ satisfying (4.5) with given constants $K \geq 8$, $0 < \rho_0 \leq R$, and*

$$u = 0 \quad \text{on} \quad (\partial_p Q) \setminus C_{\rho_0/2}(Y).$$

*Then the function*

(4.6) $$f_2(\rho) = \sup_{S_\rho}(u^-),$$

*where* $S_\rho = Q \cap (\partial_x C_{K\rho,\rho}(Y)) = Q \cap \{|x - y| = K\rho, |t - s| < \rho^2\}$, *satisfies*

(4.7) $$f_2(\rho) \leq \left(\frac{2\rho_0}{\rho}\right)^{\gamma_2} f_2(\rho_0) \quad \text{for} \quad 0 < \rho_0 \leq \rho \leq R,$$

*where the constant* $\gamma_2 = \gamma_2(n, \nu, m, K) \to \infty$ *as* $K \to \infty$. *If* $\Omega = \mathbb{R}^n$, *the constant* $\gamma_2$ *does not depend on* $m$.

*Proof.* Since $u \equiv 0$ on $(\partial_p Q) \cap \{t \leq s - \rho_0^2\}$, we also have $u \equiv 0$ on $Q \cap \{t \leq s - \rho_0^2\}$; hence without loss of generality we may assume $y = 0$, $0 \leq s \leq \rho_0^2$. By the maximum principle applied to $u$ in $Q \setminus C_{K\rho,\rho}(Y)$, the function $f_2(\rho)$ decreases on $[\rho_0, R]$, and therefore, (4.7) holds for $\rho_0 \leq \rho \leq 2\rho_0$. For arbitrary $\rho \in (2\rho_0, R]$, there exists $Z = (z, \tau) \in S_\rho$ such that $f_2(\rho) = u^-(Z)$. Since $|z| = K\rho$, we have

$$Z \in \partial_x Q_{2\rho}(Z_0), \quad \text{where} \quad Z_0 = (z, 0) \in \partial_t Q.$$

and by Corollary 4.1,

$$f_2(\rho) = u^-(Z) \leq \sup_{Q_{2\rho}(Z_0)}(u^-) \leq \theta \sup_{Q_{4\rho}(Z_0)}(u^-).$$

Notice that two sets $Q_{4\rho}(Z_0)$ and $S_{\rho_0}$ are separated by the cylindrical surface

$$S = \{|x| = (K-4)\rho\} = \{|x| = qK\rho\} \supset S_{q\rho},$$

where $q = (K-4)/K \in [1/2, 1)$, and $u \geq 0$ on $S \setminus S_{q\rho}$. By the maximum principle, we obtain

$$f_2(\rho) \leq \theta \sup_{S_{q\rho}}(u^-) = \theta f_2(q\rho) = q^{\gamma_2} f_2(q\rho),$$

where $\gamma_2 = \log_q \theta > 0$. Now we choose $k \geq 1$ satisfying $\rho_0 \leq q^k \rho \leq 2\rho_0$, and using iteration, we get the desired estimate (4.7):

$$f_2(\rho) \leq q^{k\gamma_2} f_2(q^k \rho) \leq \left(\frac{2\rho_0}{\rho}\right)^{\gamma_2} f_2(\rho_0).$$

Finally, for $K \geq 8$ we have

$$\frac{1}{q} = 1 + \frac{4}{K-4} \leq 1 + \frac{8}{K}, \quad \ln\left(\frac{1}{q}\right) \leq \frac{8}{K},$$

$$\gamma_2 = \log_q \theta = \frac{\ln(1/\theta)}{\ln(1/q)} \geq \frac{K \ln(1/\theta)}{8} \to \infty \quad \text{as} \quad K \to \infty.$$

Lemma 4.4 is proved. $\square$

4.2. *Proof of Theorem* 1.1. We may assume $r = 1$ and $K \geq 8$ is large enough to guarantee the inequality $\gamma_1 < \gamma_2$ between two constants $\gamma_1$ and $\gamma_2$ in Lemmas 4.3(a) and 4.4.

Using Lemma 4.1(a) and then the Harnack inequality, we get the estimate

$$(4.8) \qquad \omega^X(\triangle_1) \geq N^{-1} \quad \text{on} \quad U_R = \{(x,t) : |x| \leq K\sqrt{t} \leq KR\},$$

where $N = N(n, \nu, K, R) > 0$. Taking $N_0 = 2N$ with this constant $N$, we have

$$u(X) = N_0 \omega^X(\triangle_1) - \omega^X(\triangle_2) \geq 2 - 1 = 1 \quad \text{on} \quad U_R.$$

We will show that from weaker estimates

$$(4.9) \qquad u \geq 1 \quad \text{on} \quad U_8, \quad u \geq 0 \quad \text{on} \quad U_{R_0},$$

with some large constant $R_0 = R_0(n, \nu, K)$ which will be specified below, it follows that

$$(4.10) \qquad u \geq 0 \quad \text{on} \quad U_\infty;$$

i.e. the desired estimate $\omega^X(\triangle_2) \leq N_0 \omega^X(\triangle_1)$ with $N_0 = N_0(n, \nu, K)$.

Suppose (4.10) fails. Then one can choose $\rho \geq R_0/2$ such that $u \geq 0$ on $U_\rho$ and $u < 0$ at some point $X = (x, 4\rho^2)$ with $|x| \leq (2K\rho)^2$. We will use the representation of $u(X)$ through the $L$-caloric measure $\omega^X$ on the parabolic boundary of the set $Q \setminus C_{K\rho,\rho}$, where

$$Q = \mathbb{R}^n \times (0, \infty), \quad C_{K\rho,\rho} = B_{K\rho}(0) \times (-\rho^2, \rho^2).$$

By Lemmas 4.3(a) and 4.4, where $\rho_0 = 4$ (such choice of $\rho_0$ will help to extend our arguments to the proof of Theorem 1.2), we have

$$u \geq (4/\rho)^{\gamma_1} \quad \text{on} \quad B_\rho^+ = B_\rho(0) \times \{\rho^2\} \subset \partial_p(Q \setminus C_{K\rho,\rho}),$$
$$u \geq -u^- \geq -(8/\rho)^{\gamma_2} \quad \text{on} \quad S_\rho = Q \cap (\partial_x C_{K\rho,\rho}) \subset \partial_p(Q \setminus C_{K\rho,\rho}),$$

and $u \geq 0$ on the remaining part of $\partial_p(Q \setminus C_{K\rho,\rho})$. Therefore,

$$\begin{aligned}
0 > u(X) &= \int_{\partial_p(Q \setminus C_{K\rho,\rho})} u \, d\omega^X \geq \int_{B_\rho^+} u \, d\omega^X + \int_{S_\rho} u \, d\omega^X \\
&\geq \omega^X(B_\rho^+) \cdot (4/\rho)^{\gamma_1} - (8/\rho)^{\gamma_2}.
\end{aligned}$$

Similarly (see (4.8)), we also have

$$(4.11) \qquad 1 \leq N_1 \omega^X(B_\rho^+) \quad \text{on} \quad B_{2K\rho}(0) \times \{4\rho^2\},$$

where $N_1 = N_1(n, \nu, K)$. Now the previous estimate implies

$$2^{-\gamma_2} \left(\frac{\rho}{4}\right)^{\gamma_2 - \gamma_1} < N_1;$$

hence $2\rho < R_0$, if we choose $R_0 = R_0(n, \nu, K) > 0$ such that

$$2^{-\gamma_2}(R_0/8)^{\gamma_2 - \gamma_1} \geq N_1.$$

By (4.9), $u \geq 0$ at $X = (x, 4\rho^2) \in U_{2\rho} \subset U_{R_0}$. But $u(X) < 0$ by the choice of $X$. This contradiction proves (4.10), so the proof of Theorem 1.1 is complete. □

4.3. *Proof of Theorem* 1.2. First we prove the following generalization of the estimate (4.11) for Lipschitz domains $\Omega$.

LEMMA 4.5. *Let* $Y = (y, s) \in \partial_p Q$. *Given the constants* $0 < 2\rho \leq \lambda r_0$,

(4.12) $\quad \omega^X(\Omega \times \{s + \rho^2\}) \leq N_1 \omega^X(\Omega_\rho^+)$ *on* $D = \Omega_{2K\rho}(y) \times \{s + 4\rho^2\}$

*with* $N_1 = N_1(n, \nu, m, \lambda, K)$, *where*

$$\Omega_\rho^+ = \Omega_\rho^+(Y) = \Omega_\rho^{\mu\rho'}(y) \times \{s + \rho^2\}, \qquad \rho' = \min(\rho, r_0).$$

*Proof.* We set $s_0 = s + 4\rho^2$, $Y_0 = (y, s_0)$, and $r = \rho'/4$. From the inequality $2\rho \leq \lambda r_0$ it follows $2\rho \leq \lambda' \rho' = 4\lambda' r$, where $\lambda' = \max(2, \lambda)$, and then

$$D = \Omega_{2K\rho}(y) \times \{s_0\} \subset Q_{4K\lambda'r, r}(Y_0)$$
$$= \Omega_{4K\lambda'r}(y) \times (s_0 - r^2, s_0 + r^2).$$

We apply Theorem 2.4, with the constant $K' = 4K\lambda'$ instead of $K$, to the functions

$$u(X) = \omega^X(\Omega_\rho^+), \qquad v(X) = \omega^X(\Omega \times \{s + \rho^2\}) \leq 1.$$

This gives

$$\sup_D \frac{v}{u} \leq \sup_{Q_{K'r,r}(Y_0)} \frac{v}{u} \leq N \sup_{\Omega_{2r}^-(Y_0)} \frac{1}{u}.$$

By assumption (2.3), the set $\Omega_\rho^{\mu\rho'}(y)$ contains a ball of radius $\mu\rho' = 4\mu r$. Using Lemma 4.1(b) and then the Harnack inequality, we obtain

$$u \geq N^{-1} \quad \text{on} \quad \Omega_{2r}^-(Y_0) = \Omega_{2r}^{2\mu r}(y) \times \{s + 4\rho^2 - 4r^2\}.$$

Therefore, $v/u \leq N_1 = N_1(n, \nu, m, \lambda, K)$ on $D$. Lemma 4.5 is proved. □

Now we begin the *proof of Theorem* 1.2. Following the lines of the proof of Theorem 1.1, we assume $r = 1$ and choose $K \geq 8$ to guarantee the inequality $\gamma_1 < \gamma_2$, where $\gamma_1$ and $\gamma_2$ are the constants in Lemmas 4.3(b) and 4.4 correspondingly. Our goal is to show that

(4.13) $\qquad u(X) = N\omega^X(\triangle_1) - \omega^X(\triangle_2) \geq 0 \quad \text{on} \quad U'_{\lambda r_0}$

with $N = N(n, \nu, m, \lambda, K)$, where

$$U'_R = \{(x, t) \in Q : |x - y| \leq K\sqrt{t - s}, \quad 4 \leq \sqrt{t - s} \leq R\}.$$

By Lemma 4.1(b) and the Harnack inequality,

(4.14) $\qquad \omega^X(\triangle_1) \geq N_0^{-1} \quad \text{on} \quad \Omega_\rho^+ \text{ for } 0 < 2\rho \leq R \leq \lambda r_0,$

with $N_0 = N_0(n, \nu, m, \lambda, K, R)$. Using (4.14), (4.12), and the comparison principle, we have

$$N_0 N_1 \omega^X(\triangle_1) \geq N_1 \omega^X(\Omega_\rho^+) \geq \omega^X(\Omega \times \{s + \rho^2\}) \geq \omega^X(\triangle_2)$$

on the set

$$\Omega_{2K\rho}(y) \times \{s + 4\rho^2\} = U'_R \cap \{t = s + 4\rho^2\},$$

for $4 \leq 2\rho \leq R \leq \lambda r_0$. Therefore, the function $u$ in (4.13), with the constant $N = N_0 N_1$ depending on $R$, satisfies $u \geq 0$ on $U'_R$. This implies the desired estimate (4.13) if $\lambda r_0$ does not exceed a constant $R_0 = R_0(n, \nu, m, \lambda, K) > 0$, which is chosen from the relation $2^{-\gamma_2}(R_0/8)^{\gamma_2-\gamma_1} \geq N_1 = N_1(n, \nu, m, \lambda, K)$, the constant in Lemma 4.5.

Now it remains to consider the case $R_0 < \lambda r_0$. By the above arguments, there exist $N = N(n, \nu, m, \lambda, K)$ such that the function $u(X)$ in (4.13) satisfies

(4.15) $\qquad u \geq 1 \quad \text{on} \quad \Omega_\rho^+, \text{ for } 4 \leq \rho \leq 8; \qquad u \geq 0 \quad \text{on} \quad U'_{R_0}.$

We will show these properties of $u$ imply (4.13), i.e. $u \geq 0$ on $U'_{\lambda r_0}$. Suppose otherwise. Then we choose $\rho > 4$ such that $u \geq 0$ on $U'_\rho$, and $u < 0$ at some point $X \in \Omega_{2K\rho}(y) \times \{s + 4\rho^2\} \subset U'_{2\rho}$. For the $L$-caloric measure $\omega^X$ on $\partial_p(Q \setminus C_{K\rho,\rho})$, from Lemmas 4.3(b) and 4.4 it follows

$$\begin{aligned} 0 &> u(X) \geq \int_{\Omega_\rho^+} u \, d\omega^X + \int_{S_\rho} u \, d\omega^X \\ &\geq \omega^X(\Omega_\rho^+) \cdot (4/\rho)^{\gamma_1} - \omega^X(S_\rho) \cdot (8/\rho)^{\gamma_2}. \end{aligned}$$

By the comparison principle and Lemma 4.5,

$$\omega^X(S_\rho) \leq \omega^X(\Omega \times \{s + \rho^2\}) \leq N_1 \omega^X(\Omega_\rho^+).$$

The previous inequalities yield

$$2^{-\gamma_2}(\rho/4)^{\gamma_2-\gamma_1} < N_1 \leq 2^{-\gamma_2}(R_0/8)^{\gamma_2-\gamma_1};$$

hence $2\rho < R_0$, $X \in U'_{2\rho} \subset U'_{R_0}$, and $u(X) \geq 0$ by virtue of (4.15). However, by the choice of $X$, $u(X) < 0$. This contradiction proves Theorem 1.2. $\square$

University of Minnesota, Minneapolis, MN
*E-mail address*: safonov@math.umn.edu

University of Texas, Austin, TX
*E-mail address*: yyuan@math.umn.edu